\documentclass{amsproc}

\usepackage{amssymb}
\usepackage{lineno,hyperref}
\usepackage{amsfonts}
\usepackage{amsthm}
\usepackage{amsmath}
\usepackage{subfig}
\usepackage{graphicx}
\usepackage{courier}

\newtheorem{theorem}{Theorem}[section]

\theoremstyle{definition}

\theoremstyle{remark}
\newtheorem{remark}[theorem]{Remark}

\numberwithin{equation}{section}

\input{tcilatex}

\begin{document}

\title[3D nonlinear SAR imaging via convexification]{Numerical reconstruction for 3D nonlinear SAR imaging via a version of the convexification method}

\author[V. A. Khoa]{Vo Anh Khoa}
\address{Department of Mathematics, Florida A\&M University, Tallahassee, FL 32307, USA}
\email{anhkhoa.vo@famu.edu}
\thanks{This work was supported by US Army Research Laboratory and US Army Research Office grant W911NF-19-1-0044. V. A. Khoa was supported by the Faculty Research Awards Program (FRAP) at Florida A\&M University, under the project \#007633.}

\author[M. V. Klibanov]{Michael Victor Klibanov}
\address{Department of Mathematics and Statistics, University of North Carolina at Charlotte, Charlotte, NC 28223, USA}
\email{mklibanv@uncc.edu}

\author[W. G. Powell]{William Grayson Powell}
\address{Department of Mathematics, University of Wisconsin Madison, Madison, WI 53706, USA}
\email{wgpowell@wisc.edu}

\author[L. H. Nguyen]{Loc Hoang Nguyen}
\address{Department of Mathematics and Statistics, University of North Carolina at Charlotte, Charlotte, NC 28223, USA}
\email{loc.nguyen@uncc.edu}

\subjclass{78A46, 65L70, 65C20}


\keywords{SAR imaging, coefficient inverse problem,  convexification, global convergence, simulated data, delay-and-sum procedure}

\begin{abstract}
This work extends the applicability of our recent convexification-based algorithm for constructing images of the dielectric constant of buried or occluded target. We are orientated towards the detection of explosive-like targets such as antipersonnel land mines and improvised explosive devices in the non-invasive inspections of buildings. In our previous work, the method is posed in the perspective that we use multiple source locations running along a line of source to get a 2D image of the dielectric function. Mathematically, we solve a 1D coefficient inverse problem for a hyperbolic equation for each source location. Different from any conventional Born approximation-based technique for synthetic-aperture radar, this method does not need any linearization. In this paper, we attempt to verify the method using several 3D numerical tests with simulated data. We revisit the global convergence of the gradient descent method of our computational approach.
\end{abstract}

\maketitle

\section{Introduction}
Synthetic-aperture radar (SAR) imaging is a commonly used technique in
reconstructing images of surfaces of planets, and in detecting antipersonnel
land mines and improvised explosive devices; cf. e.g. \cite%
{Amin2011,Carn1999,Gilman2017,Rotheram1985} for some essential backgrounds
concerning SAR imaging. The conventional SAR imaging is
based on the linearization via the Born approximation \cite{Gilman2017}.
While providing accurate images of shapes and locations of targets, the
linearization is problematic to produce accurate values of their dielectric
constants, see, e.g. \cite{Klibanov2021}. However, since dielectric
constants characterize constituent materials of targets, then it is
obviously important to accurately compute them.

To address this problem, a new convexification based nonlinear SAR
imaging, CONSAR, was proposed in \cite{Klibanov2021} for the first time.
Testing of CONSAR for both computationally simulated and experimentally
collected data has shown that dielectric constants of targets are indeed
accurately computed \cite{Klibanov2021,Klibanov2021a}. However, these
publications are only about obtaining 2D SAR images on the so-called
\textquotedblleft slant range" plane. \ In this case both the point source
and the radiating antenna move along an interval of a single straight line.
Thus, unlike \cite{Klibanov2021,Klibanov2021a}, the current paper is about
testing the CONSAR technique for 3D SAR images. By our experience, collection of experimental SAR data along just one line is time consuming 
\cite{Klibanov2021a}. Therefore, we consider here the case when the source
and the detector concurrently move along only three (3) parallel lines. In
principle, one should have many such lines of course to get accurate
solutions. Hence, because of an obvious lack of the information content in
the input data, it is unrealistic to expect that resulting SAR images would
be accurate in terms of shapes of targets. This is why sizes of computed
images of our targets are quite different from the true ones. Nevertheless,
we accurately image dielectric constants and locations of targets.
Only computationally simulated data are considered here. 

In this paper, our target application is the non-invasive inspections of
buried targets when one tries to image the buried or occluded target using
backscattering data collected by e.g., flying devices or roving vehicles
(cf. \cite{Amin2011,Nguyen2008}). Having knowledge of the dielectric
constant of the target plays a vital role in the development of future
classification algorithms that distinguish between explosives and clutter in
military applications. Thus, it is anticipated that the knowledge of
dielectric constants of targets should decrease the false alarm rate.

In fact, the conventional time dependent SAR data when the source
and the detector move concurrently along an interval of a straight line,
depend on two variables. On the other hand, we want to image 3D objects.
Therefore, the data for the corresponding Coefficient Inverse Problem (CIP) are underdetermined. This is
why we work here with the above mentioned more informative case of three
parallel lines. Although these data are also underdetermined, we still hope
to reconstruct some approximations of 3D images of targets. It is worth
mentioning that solving underdetermined or non overdetermined CIPs is
omnipresent in applications that we have mentioned above. This is one of the
main challenges that researchers encounter in practice is finding a good
approximation of the target's dielectric constant. This is also a long
standing research project that the corresponding author and his research
team have been working on for almost a decade, see, e.g. \cite%
{Klibanov2015,Thanh2015,Kolesov2017,Nguyen2018,Klibanov2019,Khoa2020,Khoa2020a}
for mathematical models and numerical methods to identify the material of
the buried targets.

Solving for the target's dielectric constant in this scenario is a CIP for a
hyperbolic PDE. It is well-known that CIPs are both severely ill-posed and
nonlinear thus, requiring a careful designation of an appropriate numerical
method. Very recently, we have proposed in \cite{Klibanov2021,Klibanov2021a}
a new convexification-based algorithm to construct a 2D image of the
target's dielectric constant using SAR data. In fact, the approach of \cite%
{Klibanov2021,Klibanov2021a} means a nonlinear SAR imaging, as opposed to
the conventional linear one. We point out that the computational approach of 
\cite{Klibanov2021,Klibanov2021a} to the SAR data is a heuristic one.
Nevertheless, it was shown in \cite{Klibanov2021} that this approach works
well numerically for computationally simulated data. It is more encouraging
that it was demonstrated in \cite{Klibanov2021,Klibanov2021a} that this
approach works well for experimentally collected data. The same approach to
the SAR data is used in this paper. Using multiple source locations running
along an interval of a straight line, our computational approach is to solve
a number of 1D CIPs.

To solve numerically each of our CIPs, we transform first the hyperbolic
equation into a nonlocal nonlinear boundary value problem (BVP), based upon
the original idea of \cite{Smirnov2020}. An approximate solution of this BVP
is then found by minimizing a suitable weighted Tikhonov-like cost
functional, which is proved to be globally strictly convex. Its strict
convexity is based on the presence of the Carleman Weight Function (CWF) in
it. The CWF is the function which is used as the weight function in the
Carleman estimate for the principal part of the differential operator
involved in that BVP. Ultimately, the 2D image which we construct for the
case of each of our three lines can be understood as the collection
of 1D cross-sections of the so-called slant-range plane established by the
\textquotedblleft moving\textquotedblright\ source/antenna
positions. To construct the 3D image, we combine those 2D ones, see
sections \ref{sec:3} and \ref{sec:experiments} for details.

Conventional numerical methods for CIPs are based on minimizations of
various least squares cost functionals, including the Tikhonov
regularization functional, see, e.g. \cite{Chavent,Gonch1,Gonch2}. It is
well known, however, that such a functional is typically non convex and has,
therefore, multiple local minima and ravines, which, in turn plague the
minimization procedure. The convexification method was initially proposed in 
\cite{KlibIous,Klib97} with the goal to avoid that phenomenon of multiple
local minima and ravines of conventional least squares cost functionals for
CIPs. While publications \cite{KlibIous,Klib97} are purely theoretical ones,
the paper \cite{Bak} has generated many follow up numerical studies of the
convexification of this research group, see, e.g. the above cited
publications of this research group and references cited therein. We also
refer to the recently published book \cite{KL}. 

The convexification uses a significantly modified idea of the
Bukhgeim--Klibanov method \cite{BukhKlib} of 1981. However, in \cite%
{BukhKlib} and many follow up publications, only proofs of global uniqueness
and stability theorems for CIPs were studied, see, e.g. \cite%
{BK,BY,KT,Ksurvey,Yam} and references cited therein. On the other hand, the
convexification is going in the numerical direction: it is designed for
constructions of \emph{globally convergent} numerical methods for CIPs. This
is unlike the conventional locally convergent numerical methods for CIPs,
which are based on minimizations of those functionals.

As stated above, depending on the differential operator involved in the CIP,
the convexification method is based on an appropriate Carleman Weight
Function to construct a weighted Tikhonov-like cost functional. The Carleman
Weight Function plays two roles. First, it convexifies the conventional
Tikhonov functional. Second, it helps to control nonlinearities involved in
the underlying CIP.

\textbf{Definition 1.} \emph{We call a numerical method for a CIP globally
	convergent if there is a theorem which claims that this method reaches at
	least one point in a sufficiently small neighborhood of the true solution of
	that CIP without any advanced knowledge of a small neighborhood of this
	solution.} 

The convexification is a globally convergent numerical method because one
can prove the convergence of the scheme towards the true solution starting
from any point located in a given bounded set of a suitable Hilbert space
and the diameter of this set can be arbitrary large. 

As to the moving sources, we remark that in the publication \cite{Khoa2020}
a version of the convexification method was developed for the first time
both analytically and numerically for a 3D inverse medium problem for the
Helmholtz equation with the backscattering data generated by point sources
running along a straight line and at a fixed frequency. Then this result was
confirmed on experimental data in \cite{Khoa2020a,Khoa2020b}. In \cite%
{Khoa2020,Khoa2020a,Khoa2020b}, we were working with the data in the
frequency domain, whereas the SAR data are the time-dependent ones. In
addition, in the case of the SAR data, both source and detector move
concurrently along a straight line. On the other hand, in \cite%
{Khoa2020,Khoa2020a,Khoa2020b} the source moves along a part of a straight
line, whereas the detector moves independently along a part of a plane. We
also refer to works of Novikov and his group \cite{Alex,Nov} for a different
approach to CIPs with the moving source.

We also add that our recent publications \cite{Klibanov2021,Klibanov2021a}
about CONSAR were focusing only on the 2-D reconstructions. Therefore, the
main difference between these papers and the current one is that here we
obtain a 3-D image of the target dielectric constant, assuming that the data
are collected from multiple lines of sources. 

Our paper is organized as follows. In sections \ref{sec:2}, \ref{sec:3} and %
\ref{sec:4}, we revisit our setting for the CIP imaging and our
convexification-based method proposed in \cite{Klibanov2021a}. In section %
\ref{sec:experiments}, we verify the 3-D numerical performance of our method
using computationally simulated data. Finally, we present concluding remarks
in section \ref{sec:6}. 

\section{Statement of the forward problem}

\label{sec:2}

Prior to the statement of the forward problem, we discuss the setting of SAR
imaging. We focus only on the stripmap SAR processing; see e.g. \cite%
{Showmana} for details of this setup. With the non-invasive inspection of
buildings and land mines detection being in mind, we prepare to have a
radiating antenna and a receiver in the SAR device. As the antenna radiates
pulses of the time-resolved component of the electric wave field, the
receiver collects the backscattering signal. Our SAR device moves along a
straight line. In this work, we assume that it moves it along multiple
parallel lines, which, we hope, may create a 3-D image of the target. It is
worth mentioning that in our inverse setting in section \ref{sec:3} below,
we assume that the antenna and the receiver coincide, and form a point,
which moves along certain parallel lines. This assumption is a typical one
in SAR imaging; cf. \cite{Gilman2017}. Indeed, in practice the antenna and
the receiver are rather close to each other, and the distance between them
is much less than the distance between the line over which they concurrently
move and the targets of interest. However, when generating the SAR data in
computational simulations, we solve a forward problem, in which we assume
that the antenna is a disk.

We denote below $\mathbf{x}=(x,y,z)\in \mathbb{R}^{3}.$ For $z^{0}>0$, we consider the source location $\mathbf{x}%
^{0}=(x^{0},y_{j}^{0},z^{0})\in L_{j}^{\text{src}}$ for $j=\overline{1,j_{0}}
$, $j_{0}\in \mathbb{N}^{\ast }$. By choosing a certain value of $j_{0}$ and
varying $x^{0}$ in the interval $(-L,L)$ for $L>0$, we define a finite set
of lines $L_{j}^{\text{src}}$ of sources/receivers as follows: 
\begin{equation}
	L_{j}^{\text{src}}:=\left\{ \mathbf{x}^{0}=\left(
	x^{0},y_{j}^{0},z^{0}\right) \in \mathbb{R}^{3}:x^{0}\in \left( -L,L\right)
	,j=\overline{1,j_{0}}\right\} .  \label{1}
\end{equation}%
Therefore, each line of sources $L_{j}^{\text{src}}$ is parallel to the $x-$%
axis. The antenna and the receiver coincide and are located at the point $%
\mathbf{x}^{0}=\left( x^{0},y_{j}^{0},z^{0}\right) \in L_{j}^{\text{src}}.$
The number $y_{j}^{0}$ is the same for all points of $L_{j}^{\text{src}}$
and this number is changing as lines $L_{j}^{\text{src}}$ changes. The
number $z^{0}$ is the same for all lines (\ref{1}). Assuming that the
parameter $x^{0}\in \left( -L,L\right) $ in (\ref{1}), we actually assume
that the antenna/receiver point characterized by $x^{0}$ moves along each
line $L_{j}^{\text{src}},$ $j=\overline{1,j_{0}}.$

Let the number $R>0.$ Our domain of interest is 
\begin{equation*}
	\Omega =(-R,R)^{3}\subset \mathbb{R}^3.
\end{equation*}%
This domain consists of two parts. The lower part $\Omega \cap \left\{
z<0\right\} $ mimics the sandbox containing the unknown buried object.
Meanwhile, the upper part $\Omega \cap \left\{ z>0\right\} $ is assumed to
be the air part, where we move our SAR device. For simplicity, we assume
that the interface between the air and sand parts is located at $\left\{
z=0\right\} $. Let $\varepsilon _{r}(\mathbf{x})\in C^{1}\left( \mathbb{R}%
^{3}\right) $ be a function that represents the spatially distributed
dielectric constant. We assume that%
\begin{equation}
	\left\{ 
	\begin{array}{c}
		\varepsilon _{r}\left( \mathbf{x}\right) \geq 4,\text{ for }\left\{
		z<0\right\} \cap \Omega , \\ 
		\varepsilon _{r}\left( \mathbf{x}\right) =1,\text{ for }\left\{ z>0\right\}
		\cap \Omega .\text{ }%
	\end{array}%
	\right.   \label{eps}
\end{equation}

Physically, assumption (\ref{eps}) is reasonable. Cf. e.g. \cite{Khoa2020a},
we know that the dielectric constant of the air/vacuum is identically one.
Meanwhile, for the other materials (for example, the dry sand which follows
the non-invasive inspections of our interest), the dielectric constant is
larger than the unity.

\begin{figure}[htbp]
	\centering
	\includegraphics[scale = 0.47]{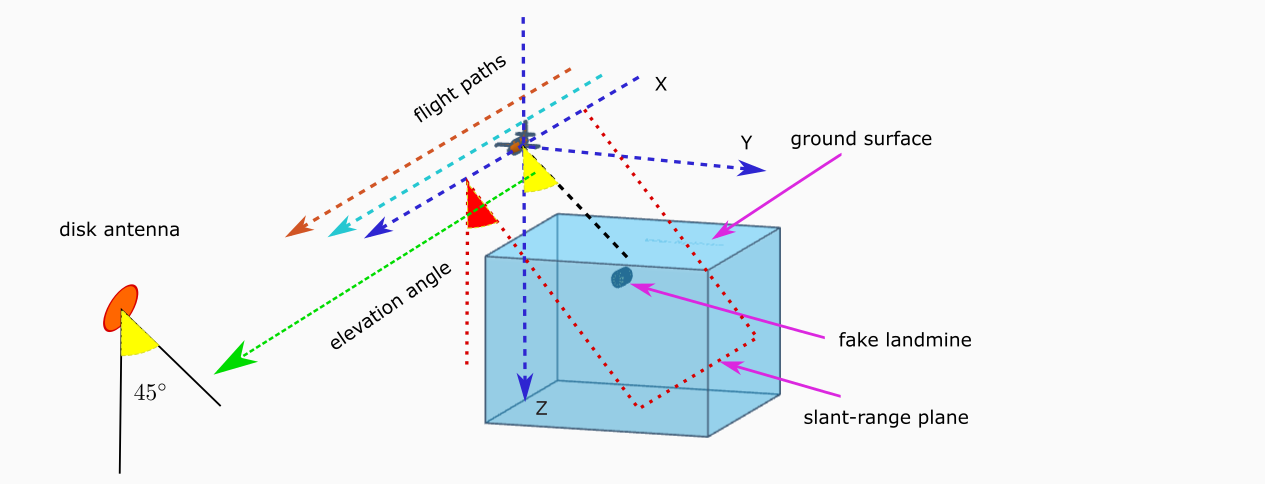}
	\caption{A schematic diagram of our SAR imaging setting. A flying device
		moves along multiple parallel straight lines in $x$-direction to collect the
		backscattering data. On this device, we assume that the transmitter and the
		detector have the same position. The transmitter is a disk-shaped antenna
		with an elevation angle of $\protect\pi/4$ (yellow part). For our purpose of
		mimicking the landmine detection, we assume the unknown object is buried in
		a dry sand region. The notion of our convexification-based method is solving
		1D CIPs for every source location. Under a certain elevation angle (red
		part), solutions of these 1-D CIPs allow us to approximately reconstruct a
		3-D image. }
	\label{fig:1}
\end{figure}

Let $T>0$ be the observation time. Practically, $T$ is found by doubling the
distance from the line of sources to the furthest point of the sandbox. For
every source location $\mathbf{x}^{0}\in L_{j}^{\text{src}}$ (see (\ref{1}%
)), we consider the following forward problem:%
\begin{align}
	\begin{cases} \label{forward1}
		\nu^{-2}\left(\mathbf{x}\right)u_{tt}=\Delta u+F\left(\mathbf{x},\mathbf{x}^{0},t\right),\\
		u\left(\mathbf{x},\mathbf{x}^{0},0\right)=u_{t}\left(\mathbf{x},\mathbf{x}^{0},0\right),\\
		\left.\left[\partial_{n}u\left(\mathbf{x},\mathbf{x}^{0},t\right)+\nu^{-1}\left(\mathbf{x}\right)u\left(\mathbf{x},\mathbf{x}^{0},t\right)\right]\right|_{\partial\Omega}=0,\\
		\mathbf{x}^{0}\in L_{j}^{\text{src}},\mathbf{x}\in\Omega,t\in\left(0,T\right),j=\overline{1,j_{0}},
	\end{cases}
\end{align}
where $n=n\left( \mathbf{x}\right) ,\mathbf{x}\in \partial \Omega $ is the
unit outward looking normal vector on $\partial \Omega .$ In (\ref{forward1}%
), $u(\mathbf{x},\mathbf{x}^{0},t)$ is the amplitude of a time-resolved
component of the electric wave field. The function $\nu (\mathbf{x})=c_{0}/%
\sqrt{\varepsilon _{r}(\mathbf{x})}$ represents the speed of light in the
medium with $c_{0}$ is the speed of light in the vacuum. Here, the source
function $F$ is defined as 
\begin{equation}
	F\left( \mathbf{x},\mathbf{x}^{0},t\right) =%
	\begin{cases}
		\text{Re}\left( \chi _{D}\left( \mathbf{x},\mathbf{x}^{0}\right) e^{-\text{i}%
			\omega _{0}t}e^{-\text{i}\alpha _{0}\left( t-\tau _{0}/2\right) ^{2}}\right) 
		& \text{if }t\in \left( 0,\tau _{0}\right] ,\mathbf{x}^{0}\in L_{j}^{\text{%
				src}}, \\ 
		0 & \text{if }t>\tau _{0},%
	\end{cases}
	\label{source}
\end{equation}%
where $\text{i}=\sqrt{-1}$. This source function models well the linear
modulated pulse, where $\omega _{0}$ is the carrier frequency, $\alpha _{0}$
is the chirp rate \cite{Gilman2017}, and $\tau _{0}$ is its duration. The
function $\chi _{D}$ in \eqref{source} represents the disk-shaped antenna of
our interest with its radius $D>0$, thickness $\widetilde{D}>0$ and centered
at a point $\mathbf{x}\in L_{j}^{\text{src}}$. We rotate the disk via a
coordinate transformation $\mathbf{x}\mapsto \mathbf{x}^{\prime }=(x^{\prime
},y^{\prime },z^{\prime })$ defined by 
\begin{equation}
	\begin{array}{rcll}
		x^{\prime } & = & x-x^{0}, &  \\ 
		y^{\prime } & = & \cos (\theta )(y-y_{j}^{0})-\sin (\theta )(z-z^{0}), &  \\ 
		z^{\prime } & = & \sin (\theta )(y-y_{j}^{0})+\cos (\theta )(z-z^{0}), &  \\ 
		&  &  & 
	\end{array}
	\label{2000}
\end{equation}%
where $\theta $ is the elevation angle. Hereby, the function $\chi _{D}$ is
given by 
\begin{equation}
	\chi _{D}(\mathbf{x},\mathbf{x}^{0})=%
	\begin{cases}
		1 & \text{if }\sqrt{{x^{\prime }}^{2}+{y^{\prime }}^{2}}<D\text{ and }%
		|z^{\prime }|\leq \widetilde{D}, \\ 
		0 & \text{otherwise.}%
	\end{cases}
	\label{3000}
\end{equation}

We consider the so-called absorbing boundary conditions in (\ref{forward1}).
Mathematically, the solution to the acoustic wave equation is sought on the
whole space. As introduced above, the computational domain $\Omega $ is,
however, finite and the boundary condition should be equipped with the
underlying PDE. Moreover, at this point, numerical simulations show the
non-physical back reflections of the waves incident at the boundaries.
Therefore, the absorbing boundary conditions are used to minimize these
spurious reflections; cf. e.g. \cite{Renaut1992} and references cited
therein.

\section{Statement of the inverse problem}

\label{sec:3}

Our application being in mind is to identify the material of the buried
object. Therefore, our inverse problem aims to compute the spatial
distribution of the dielectric constant $\varepsilon _{r}(\mathbf{x})$
involved in the coefficient $\nu (\mathbf{x})$ in the PDE of \eqref{forward1}%
.

\subsection{Statement of the problem}

\label{sec:3.1}

\textbf{SAR Inverse Problem.} \emph{Suppose that functions }$G_{j}(\mathbf{x}%
^{0},t)=u(\mathbf{x}^{0},\mathbf{x}^{0},t)$\emph{\ for }$t\in \left(
0,T\right) ,$\emph{\ }$\mathbf{x}^{0}\in L_{j}^{\text{scr}}$\emph{\ and }$j=%
\overline{1,j_{0}}$\emph{\ are given as our SAR data. Find a function }$%
\varepsilon _{r}(\mathbf{x})\in C^{1}\left( \mathbb{R}^{3}\right) $\emph{\
	satisfying conditions in (\ref{eps}). }

Cf. \cite{Klibanov2021,Klibanov2021a}, our strategy is to image the unknown
function on the so-called slant-range plane, denoted by $P$. We define this
plane in the following manner. For each antenna/source location $\mathbf{x}%
_{j}^{0,m}=\left( x_{m}^{0},y_{j}^{0},z^{0}\right) \in L_{j}^{\text{scr}}$,
we define the central line of the antenna $CRL(\mathbf{x}_{j}^{0,m})$ that
passes through the point $\mathbf{x}_{j}^{0,m}$, $m=1,...,N,$ where $N$ is
the total number of locations of the source on the line $L_{j}^{\text{scr}}.$
Each line $CRL(\mathbf{x}_{j}^{0,m})$ is orthogonal to the $x$-axis and has
a certain angle with the plane $\left\{ y=y_{j}^{0}\right\} $. Recall that
the line $L_{j}^{\text{scr}}$ is parallel to the $x-$axis, see (\ref{1}).
Denoted by $\theta $, that angle is called the elevation angle presumably
determined by the propagating direction of the radiated energy of the
antenna, see (\ref{source})-(\ref{3000}). Henceforth, for each line $L_{j}^{%
	\text{scr}},$ the slant-range plane $P_{j}$ is defined as the plane passing
through the central lines of antennas and the line of sources moving along $%
L_{j}^{\text{scr}}$. Hence, multiple lines of sources give multiple
slant-range planes. As a by-product, the solution obtained on this plane is
called the slant-range image $\widetilde{\varepsilon }_{r,j}(\mathbf{x})$ of
the unknown dielectric constant $\varepsilon _{r}(\mathbf{x})$ on the
slant-range plane $P_{j}$. To obtain an approximate image of the function $%
\varepsilon _{r}(\mathbf{x}),$ we then combine slant-range images $%
\varepsilon _{r,j}(\mathbf{x})$ on planes $P_{j}$, $j=\overline{1,j_{0}}.$
For clarity, we depict an illustration of our SAR imaging setting in Figure %
\ref{fig:1}.

To obtain the image in the slant range plane $P_{j},$ we solve many 1-D
Coefficient Inverse Problems (CIPs) along lines of antennas $CRL(\mathbf{x}%
_{j}^{0,m}),m=1,...,N.$ Let $\widetilde{\varepsilon }_{r,j,m}(\mathbf{x})$
be the solution of this CIP along the line $CRL(\mathbf{x}_{j}^{0,m}).$ Then
the resulting 2-D function $\widetilde{\varepsilon }_{r,j}(\mathbf{x})$ is
obtained by averaging of these functions over $m=1,...,N,$ see details in
section \ref{sec:experiments}.

\subsection{1D Coefficient Inverse Problems}

\label{sec:3.2}

We now present our computational approach to calculate the slant-range
dielectric constant $\widetilde{\varepsilon }_{r,j}(\mathbf{x})$ on each
slant range plane $P_{j}$. For every $j=\overline{1,j_{0}}$, we consider the 
$m-$th source location $\mathbf{x}_{j}^{0,m}$ for $m=\overline{1,N}$.
Observe that when introducing a pseudo variable $x\in \mathbb{R}$ and
considering a smooth function $c(x)$ for our unknown function $\widetilde{%
	\varepsilon }_{r}(\mathbf{x})$ along each central line $CRL(\mathbf{x}%
_{j}^{0,m})$ of the antenna passing through the source $\mathbf{x}_{j}^{0,m}$%
, we can scale to $x\in (0,1)$ based upon the length of that central line.
This way, we treat $c(x)$ as the solution of a 1-D CIP solved along each
central line of the source $\mathbf{x}_{j}^{0,m}$.

By (\ref{eps}), for a number $\overline{c}\geq 4$ we consider the function $%
c\left( x\right) \in C^{3}\left( \mathbb{R}\right) $ such that 
\begin{equation}
	c\left( x\right) =%
	\begin{cases}
		\in \left[ 4,\overline{c}\right]  & \text{for }x\in \left( 0,1\right) , \\ 
		1 & \text{for }x\leq 0\text{ and }x\geq 1.%
	\end{cases}
	\label{99}
\end{equation}%
Hereby, we consider the following forward problem: 
\begin{equation}
	\begin{cases}
		\frac{c(x)}{c_{0}^{2}}u_{tt}=u_{xx} & \text{for }x\in \mathbb{R},t\in \left(
		0,T\right) , \\ 
		u\left( x,0\right) =0,\quad u_{t}\left( x,0\right) =\delta \left( x\right) 
		& \text{for }x\in \mathbb{R}.%
	\end{cases}
	\label{forward}
\end{equation}%
In (\ref{forward}), the variable $t$ is considered in nanosecond (ns),
whilst the unit for $x$ is meter (m). Physically, the speed of light in the
vacuum $c_{0}=0.3$ (m/ns). Considering dimensionless variables $X$ and $\tau
,$ 
\begin{equation}
	X=x/(0.3\text{m})\text{ and }\tau =t/\text{ns},  \label{4000}
\end{equation}
we arrive at 
\begin{equation}
	u_{XX}=(0.3\text{m})^{2}u_{xx},\quad u_{\tau \tau }=(1\text{ns})^{2}u_{tt}.
	\label{trans}
\end{equation}%
Thus, we obtain the dimensionless regime of the forward problem (\ref%
{forward}) becomes:%
\begin{equation}
	\left\{ 
	\begin{array}{c}
		\widetilde{c}(X)u_{\tau \tau }=u_{XX},\quad \widetilde{c}(X)=c(0.3X), \\ 
		u\left( X,0\right) =0,u_{\tau }\left( X,0\right) =\delta \left( X\right) .%
	\end{array}%
	\right.   \label{100}
\end{equation}%
The change of variables (\ref{4000}) means that $X=1$ implies $x=0.3$ (m) in
the reality, and $\tau =1$ indicates $t=1$ (ns). We now state our
coefficient inverse problem.

\textbf{1-D Coefficient Inverse Problem (1-D CIP)}. \emph{Assume that the
	following functions }$g_{0}\left( \tau \right) $\emph{\ and }$g_{1}\left(
\tau \right) \ $\emph{\ are given: }%
\begin{equation}
	u\left( 0,\tau \right) =g_{0}\left( \tau \right) ,\quad u_{X}\left( 0,\tau
	\right) =g_{1}\left( \tau \right) ,\quad \tau \in \left( 0,T\right) .
	\label{3.4}
\end{equation}%
\emph{Determine the function }$\widetilde{c}\left( X\right) $\emph{\ in %
	\eqref{100} such that the corresponding function }$c\left( x\right) \in
C^{3}\left( \mathbb{R}\right) $\emph{\ satisfies conditions (\ref{99}).}

To solve this 1D CIP numerically, we employ a version of \cite%
{Klibanov2021,Klibanov2021a} of the convexification method. First, we change
variables as: 
\begin{equation}
	Y=Y(X)=\int_{0}^{X}\sqrt{\widetilde{c}(s)}ds.  \label{tau}
\end{equation}%
Observe that $dY/dX=\sqrt{\widetilde{c}(X)}\geq 1$. Consider the following
functions: 
\begin{align}
	& w\left( Y,\tau \right) =u\left( X\left( Y\right) ,\tau \right) \widetilde{c}%
	^{1/4}\left( X\left( Y\right) \right) ,\quad Q\left( Y\right) =\widetilde{c}%
	^{-1/4}\left( X\left( Y\right) \right) ,\nonumber \\ & p\left( Y\right) =\frac{%
		Q^{\prime \prime }\left( Y\right) }{Q\left( Y\right) }-2\left[ \frac{%
		Q^{\prime }\left( Y\right) }{Q\left( Y\right) }\right] ^{2}.  \label{Qp}
\end{align}%
The function $p(Y)$ is smooth and $p(Y)=0$ for $Y<0$ and $Y>\sqrt{\overline{c%
}}$. Moreover, the following PDE holds with initial and boundary conditions
holds: 
\begin{align}
	& w_{\tau \tau }=w_{YY}+p\left( Y\right) w\quad \text{for }Y\in \mathbb{R}%
	,\tau \in \left( 0,\widetilde{T}\right) ,  \label{w1} \\
	& w\left( Y,0\right) =0,\quad w_{\tau }\left( Y,0\right) =\delta \left(
	Y\right) \quad \text{for }Y\in \mathbb{R},  \label{w2} \\
	& w\left( 0,\tau \right) =g_{0}\left( \tau \right) ,\quad w_{Y}\left( 0,\tau
	\right) =g_{1}\left( \tau \right) \quad \text{for }\tau \in \left( 0,%
	\widetilde{T}\right) ,  \label{w3}
\end{align}%
where the number $\widetilde{T}\geq 2\sqrt{\overline{c}}$ depends on $T$.

\begin{remark}
	\label{rem4}
	
	Cf. \cite{Klibanov2021,Klibanov2021a} in this SAR inception, the smoothness
	of the unknown dielectric constant $c(x)$ should be $C^{3}(\mathbb{R})$.
	This assumption is used when we revisit our theorems below. In this regard,
	the function $p(Y)$ belongs to $C^{1}(\mathbb{R})$, see (\ref{Qp}). The
	number $\widetilde{T}\geq 2\sqrt{\overline{c}}$ is conditioned because it
	guarantees the uniqueness of the CIP (\ref{w1})--(\ref{w3}), see \cite[%
	Theorem 2.6 of Chapter 2]{Romanov2019}. As soon as the function $p(Y)$ is
	reconstructed, the original function $c(x)$ is reconstructed using (\ref{tau}%
	) and (\ref{Qp}) via the procedure described in \cite[section 7.2]%
	{Smirnov2020}.
\end{remark}

From now onward, we follow a novel transformation commenced in \cite%
{Smirnov2020} to obtain a nonlocal nonlinear PDE for (\ref{w1})--(\ref{w3}).
When doing so, we take an arbitrary $b\geq \sqrt{\overline{c}}$ and then
take $\mu \in (0,2\alpha b)$ for $\alpha \in (0,1/2)$. We introduce the
following 2D regions: 
\begin{align}
	& \mathcal{K}:=\left\{ \left( Y,\tau \right) \in \mathbb{R}^{2}:Y\in \left(
	0,b\right) ,\tau \in \left( 0,\widetilde{T}\right) \right\} ,  \label{R1} \\
	& \mathcal{K}_{\alpha ,\mu ,b}:=\left\{ \left( Y,\tau \right) \in \mathbb{R}%
	^{2}:Y+\alpha \tau <2\alpha b-\mu \;\text{for }Y,\tau >0\right\} .
	\label{R2}
\end{align}%
Henceforth, we consider a new function $v(Y,\tau ),$ 
\begin{equation}
	v(Y,\tau )=w(Y,\tau +Y)\in C^{3}(Y\geq 0,\tau \geq 0).  \label{200}
\end{equation}
Using multi-variable chain rules, we obtain the following PDE for $v(Y,\tau
) $ from (\ref{w1}): 
\begin{equation}  \label{v1}
	v_{YY}-2v_{Y\tau }+p\left( Y\right) v=0\quad \text{for }\left( Y,\tau
	\right) \in \mathcal{K}.
\end{equation}%
Furthermore, by (\ref{w3}), it yields $v\left( Y,0\right) =1/2$ for $Y\in
\left( 0,b\right) $, which then results in the fact that 
\begin{equation}  \label{pp}
	p(Y)=4v_{Y\tau }(Y,0),\quad Y\in (0,b).
\end{equation}

Differentiating both sides of (\ref{v1}) with respect to $\tau $ and setting 
\begin{equation}
	V(Y,\tau )=v_{\tau }(Y,\tau )\in C^{2}(\overline{\mathcal{K}}),  \label{201}
\end{equation}%
we obtain the following overdetermined boundary value problem for a
nonlinear hyperbolic equation: 
\begin{align}
	& V_{YY}-2V_{Y\tau }+4V_{Y}\left( Y,0\right) V=0\quad \text{for }\left(
	Y,\tau \right) \in \mathcal{K},  \label{V1} \\
	& V\left( 0,\tau \right) =q_{0}\left( \tau \right) :=g_{0}^{\prime }\left(
	\tau \right) \quad \text{for }\tau \in \left( 0,\widetilde{T}\right) ,
	\label{V2} \\
	& V_{Y}\left( 0,\tau \right) =q_{1}\left( \tau \right) :=g_{0}^{\prime
		\prime }\left( \tau \right) +g_{1}^{\prime }\left( \tau \right) \quad \text{%
		for }\tau \in \left( 0,\widetilde{T}\right) ,  \label{V3} \\
	& V_{Y}\left( b,\tau \right) =0\quad \text{for }\tau \in \left( 0,\widetilde{%
		T}\right) .  \label{V4}
\end{align}%
Boundary condition (\ref{V4}) follows from the absorbing boundary condition,
which was proven in \cite{Smirnov2020a} for our 1-D case.

\begin{remark}
	Essentially, our heuristic computational approach allows us to solve the
	non-overdetermined 3D CIP via the solutions of many overdetermined boundary
	value problems. The transformations we have presented are necessary because
	we want to use a suitable Carleman estimate. This estimate is applied to
	deal with nonlinearities involved in the PDE operator. Frequently, the
	nonlinearity is terms are not avoidable when working on underdetermined. and
	non-overdetermined CIPs (cf. e.g. \cite{Klibanov2015,Klibanov2019,Khoa2020}).
	
	In our SAR perspective, we observe that after computing $V(Y,\tau )$ for $%
	(Y,\tau )\in \mathcal{K}$ from (\ref{V1})--(\ref{V4}), we then obtain 
	\begin{equation}
		p(Y)=4V_{Y}(Y,0),\quad Y\in (0,b),  \label{pp2}
	\end{equation}%
	by (\ref{pp}). Hence, the target unknown coefficient $c(x)$ of the original
	1D CIP is sought after getting $Q(Y)$ in (\ref{Qp}).
\end{remark}

\section{A version of the convexification method: theorems revisited}

\label{sec:4}

Consider the subspace $H_{0}^{k}(\mathcal{K})\subset H^{k}(\mathcal{K})$ for 
$k\in \mathbb{N}^{\ast }$. In this work, we introduce 
\begin{align*}
	& H_{0}^{2}\left( \mathcal{K}\right) :=\left\{ u\in H^{2}\left( \mathcal{K}%
	\right) :u\left( 0,\tau \right) =u_{Y}\left( 0,\tau \right) ,u_{Y}\left(
	b,\tau \right) =0\right\} , \\
	& H_{0}^{4}\left( \mathcal{K}\right) :=H^{4}\left( \mathcal{K}\right) \cap
	H_{0}^{2}\left( \mathcal{K}\right) .
\end{align*}%
As mentioned above, several transformations lead us to an overdetermined
BVP. Since this problem is nonlinear, then any conventional
optimization-based approach, including the Tikhonov regularization and the
least-square method may suffer from the phenomenon of local minima and
ravines. This is the main reason which has originally prompted the
corresponding author to develop the convexification approach \cite%
{KlibIous,Klib97}. The main ingredient of the method and its variants is
using an appropriate Carleman Weight Function to \textquotedblleft
convexify\textquotedblright\ the cost functional. The associated high-order
regularization terms usually play the role in controlling nonlinear terms
involved in the PDE operator. Thereby, one obtains a unique minimizer and
then, the global convergence of the gradient-like minimization procedure is
attained. Depending on the principal parts of differential operators, one
can choose different Carleman Weight Functions, see, e.g. \cite{KL} for a
variety of Carleman estimates for partial differential operators. In this
work, we consider the following Carleman Weight Function for the linear
operator $\partial _{Y}^{2}-2\partial _{Y}\partial _{\tau }$ (see (\ref{V1}%
)): 
\begin{equation*}
	\psi _{\lambda }(Y,\tau )=e^{-2\lambda (Y+\alpha \tau )},\quad \lambda \geq
	1,\alpha \in (0,1/2).
\end{equation*}%
Denoting 
\begin{equation}
	\mathcal{S}(V)=V_{YY}-2V_{Y\tau }+4V_{Y}(Y,0)V,\text{ }(Y,\tau )\in \mathcal{%
		K},  \label{202}
\end{equation}%
we come up with the following weighted Tikhonov-like functional $J_{\lambda
	,\gamma }:H^{4}(\mathcal{K})\rightarrow \mathbb{R}_{+}$: 
\begin{equation}
	J_{\lambda ,\gamma }\left( V\right) =\int_{\mathcal{K}}\left[ \mathcal{S}%
	\left( V\right) \right] ^{2}\psi _{\lambda }dYd\tau +\gamma \left\Vert
	V\right\Vert _{H^{4}\left( \mathcal{K}\right) }^{2},\quad \gamma \in (0,1).
	\label{JJ}
\end{equation}

We choose the high-order regularization term in $H^{4}(\mathcal{K})$ because
of the embedding $H^{4}(\mathcal{K})\subset C^{2}(\overline{\mathcal{K}})$.
This $C^{2}$ regularity is helpful in controlling our nonlinear term in (\ref%
{V1}); see \cite[Theorem 2]{Smirnov2020}. Following our convexification
framework in, e.g., \cite{Khoa2020,Klibanov2019,Smirnov2020}, we introduce
the set $B\left( r,q_{0},q_{1}\right) ,$ 
\begin{align}
	& B\left( r,q_{0},q_{1}\right)   \label{ball} \\
	& =\left\{ u\in H^{4}\left( \mathcal{K}\right) :u\left( 0,\tau \right)
	=q_{0}\left( \tau \right) ,u_{Y}\left( 0,\tau \right) =q_{1}\left( \tau
	\right) ,u_{Y}\left( b,\tau \right) =0,\left\Vert u\right\Vert _{H^{4}\left( 
		\mathcal{K}\right) }<r\right\} .  \notag
\end{align}%
in which $r>0$ is an arbitrary number. However, in our computational
practice we use the $H^{2}\left( \mathcal{K}\right) -$norm, see section 5.1.

\textbf{Minimization Problem.} \emph{Minimize the cost functional }$%
J_{\lambda ,\gamma }\left( V\right) $\emph{\ defined in (\ref{JJ}) on the
	set }$\overline{B\left( r,q_{0},q_{1}\right) }.$

We are now in the position to go over our analysis in \cite%
{Smirnov2020,Smirnov2020a} of the convexification-based method for the
nonlinear boundary value problem (\ref{V1})--(\ref{V4}).

\begin{theorem}[Carleman estimate \protect\cite{Smirnov2020a}]
	\label{carleman} There exist constants $C=C(\alpha ,\mathcal{K})>0$ and $%
	\lambda _{0}=\lambda _{0}(\alpha ,\mathcal{K})\geq 1$ such that for all
	functions $u\in H_{0}^{2}(\mathcal{K})$ and all $\lambda \geq \lambda _{0}$
	the following Carleman estimate holds true: 
	\begin{align*}
		& \int\limits_{\mathcal{K}}\left( u_{YY}-2u_{Y\tau }\right) ^{2}\psi
		_{\lambda }dYd\tau \geq C\int\limits_{\mathcal{K}}\left( \lambda \left(
		u_{Y}^{2}+u_{\tau }^{2}\right) +\lambda ^{3}u^{2}\right) \psi _{\lambda
		}dYd\tau  \\
		& +C\int\limits_{0}^{b}\left( \lambda u_{Y}^{2}+\lambda ^{3}u^{2}\right)
		\left( Y,0\right) e^{-2\lambda Y}dY-Ce^{-2\lambda \widetilde{T}%
		}\int\limits_{0}^{b}\left( \lambda u_{Y}^{2}+\lambda ^{3}u^{2}\right) \left(
		Y,\widetilde{T}\right) dY.
	\end{align*}
\end{theorem}

\begin{theorem}[Global strict convexity \protect\cite{Smirnov2020}]
	\label{sconvex} For any $\lambda ,\gamma >0$ and functions $V\in \overline{%
		B\left( r,q_{0},q_{1}\right) }$, the cost functional $J_{\lambda ,\gamma }$
	defined in (\ref{JJ}) has the Fr\'{e}chet derivative $J_{\lambda ,\gamma
	}^{\prime }\left( V\right) \in H_{0}^{4}\left( \mathcal{K}\right) .$ Let $%
	\lambda _{0}=\lambda _{0}(\alpha ,\mathcal{K})\geq 1$ be the constant of
	Theorem \ref{carleman} and let $\widetilde{T}\geq 2b.$ Then there exists a
	sufficiently large number $\lambda _{1}=\lambda _{1}\left( \mathcal{K}%
	,r,\alpha \right) \geq \lambda _{0}$ such that for all $\lambda \geq \lambda
	_{1}$ and all $\gamma \in \lbrack 2e^{-\lambda \alpha \widetilde{T}},1),$
	the cost functional (\ref{JJ}) is strictly convex on the set $\overline{%
		B\left( r,q_{0},q_{1}\right) }$. Moreover, for all $V_{1},V_{2}\in \overline{%
		B\left( r,q_{0},q_{1}\right) }$ the following estimate holds true: 
	\begin{align}
		& J_{\lambda ,\gamma }\left( V_{2}\right) -J_{\lambda ,\gamma }\left(
		V_{1}\right) -J_{\lambda ,\gamma }^{\prime }\left( V_{1}\right) \left(
		V_{2}-V_{1}\right) \geq Ce^{-2\lambda \left( 2\alpha b-\mu \right)
		}\left\Vert V_{2}-V_{1}\right\Vert _{H^{1}\left( \mathcal{K}_{\alpha ,\mu
				,b}\right) }^{2} \\
		& +Ce^{-2\lambda \left( 2\alpha b-\mu \right) }\left\Vert V_{2}\left(
		y,0\right) -V_{1}\left( y,0\right) \right\Vert _{H^{1}\left( 0,2\alpha b-\mu
			\right) }^{2}+\frac{\gamma }{2}\left\Vert V_{2}-V_{1}\right\Vert
		_{H^{4}\left( \mathcal{K}\right) }^{2},  \notag
	\end{align}%
	where the constant $C=C\left( \alpha ,\mu ,b,r\right) >0$ depends only on
	listed parameters. Furthermore, the functional $J_{\lambda ,\gamma }\left(
	V\right) $ has unique minimizer $V_{\min }$ on the set $\overline{B\left(
		r,q_{0},q_{1}\right) }$ and 
	\begin{equation*}
		J_{\lambda ,\gamma }^{\prime }\left( V_{\min }\right) \left( V-V_{\min
		}\right) \geq 0,\text{ }\forall V\in \overline{B\left( r,q_{0},q_{1}\right) }%
		.
	\end{equation*}
\end{theorem}

By the conventional regularization theory (cf. \cite{Tikhonov1995}), the
existence of the ideal noiseless data $g_{0}^{\ast }\left( \tau \right)
,g_{1}^{\ast }\left( \tau \right) $ is assumed a priori. The existence of
the corresponding exact coefficient $p^{\ast }\left( Y\right) $ and the
exact function $w^{\ast }\left( Y,\tau \right) $ in (\ref{w1})--(\ref{w3})
is assumed as well. Moreover, this function $p^{\ast }\left( Y\right) $
should also satisfy conditions $p^{\ast }(Y)=0$ for $Y<0$ and $Y>\sqrt{%
	\overline{c}}$ as ones for $p(Y)$. Having the function $w^{\ast }\left(
Y,\tau \right) $, we apply the above transformations (\ref{200}), (\ref{201}%
) to obtain the corresponding function $V^{\ast }\left( Y,\tau \right) $.
Henceforth, it is natural to assume that functions $q_{0}^{\ast }$ and $%
q_{1}^{\ast }$ are the noiseless data $q_{0}$ and $q_{1}$ respectively; see (%
\ref{V2}) and (\ref{V3}). Thus, we assume below that the function $V^{\ast
}\in B(r,q_{0}^{\ast },q_{1}^{\ast })$. Hence, it follows from (\ref{V1})--(%
\ref{V4}) and (\ref{202}) that 
\begin{equation*}
	\mathcal{S}(V^{\ast })=0,\quad (Y,\tau )\in \mathcal{K},\;\text{and }p^{\ast
	}(Y)=4v_{Y}^{\ast }(Y,0),\quad Y\in (0,b).
\end{equation*}

Consider a sufficiently small number $\sigma \in (0,1)$ that characterizes
the noise level between the data $q_{0},q_{1}$ and $q_{0}^{\ast
},q_{1}^{\ast }$. To obtain the zero boundary conditions at $\left\{
x=0\right\} $ in (\ref{V2}) and (\ref{V3}), consider two functions $F\in
B(r,q_{0},q_{1})$ and $F^{\ast }\in B(r,q_{0}^{\ast },q_{1}^{\ast })$. We
assume that 
\begin{equation}
	\left\Vert F-F^{\ast }\right\Vert _{H^{4}(\mathcal{K})}<\sigma .
	\label{3.23}
\end{equation}

Consider functions $W,W^{\ast },$ 
\begin{equation*}
	W^{\ast }=V^{\ast }-F^{\ast },W=V-V\in B_{0}(2r):=\left\{ u\in H_{0}^{4}(%
	\mathcal{K}):\left\Vert u\right\Vert _{H^{4}\left( \mathcal{K}\right)
	}<2r\right\} .
\end{equation*}%
Besides, the function $W+F\in B(3r,q_{0},q_{1})$, $\forall W\in B_{0}(2r)$.
We modify the cost functional $J_{\lambda ,\gamma }\left( V\right) $ as 
\begin{equation*}
	\widetilde{J}_{\lambda ,\gamma }(W)=J_{\lambda ,\gamma }(W+F),\text{ }%
	\forall W\in B_{0}(2r).
\end{equation*}

It follows from Theorem 2 that the functional $\widetilde{J}_{\lambda
	,\gamma }$ is also globally strict convex on the ball $\overline{B_{0}(2r)}$
for $\lambda \geq \lambda _{2}=\lambda _{1}\left( \mathcal{K},3r,\alpha
\right) \geq \lambda _{1}\left( \mathcal{K},r,\alpha \right) .$ Also, by
Theorem 2, for each value of the parameter $\lambda \geq \lambda _{2},$ the
functional $\widetilde{J}_{\lambda ,\gamma }\left( W\right) $ has a unique
minimizer $W_{\min ,\lambda ,\gamma }$ on the set $\overline{B_{0}(2r)}$ and 
\begin{equation}
	\widetilde{J}_{\lambda ,\gamma }^{\prime }(W_{\min ,\lambda ,\gamma })\left(
	W-W_{\min ,\lambda ,\gamma }\right) \geq 0,\text{ }\forall W\in \overline{%
		B_{0}(2r)}.  \label{203}
\end{equation}%
It follows from the proof of Theorem 5 of \cite{Smirnov2020} that inequality
(\ref{203}) plays an important role in the proof of Theorem 3.

\begin{theorem}[Accuracy estimate of the minimizer \protect\cite{Smirnov2020}%
	]
	\label{thm:10} Assume that (\ref{3.23}) holds true and let $\widetilde{T}%
	\geq 4b.$ We choose 
	\begin{equation*}
		\beta =\frac{\alpha \left( \widetilde{T}-4b\right) +\mu }{2\left( 2\alpha
			b-\mu \right) },\quad \rho =\frac{1}{2}\min \left\{ \beta ,1\right\} .
	\end{equation*}%
	Let $\lambda _{1}=\lambda _{1}\left( \mathcal{K},r,\alpha \right) $ by the
	number of obtained Theorem \ref{sconvex} and let $\lambda _{2}=\lambda
	_{1}\left( \mathcal{K},3r,\alpha \right) \geq \lambda _{1}\left( \mathcal{K}%
	,r,\alpha \right) $. Let $\sigma _{0}\in \left( 0,1\right) $ be a
	sufficiently small number such that $\ln \sigma _{0}^{-1/\left( 2\left(
		2\alpha b-\mu \right) \right) }\geq \lambda _{2}$. For any $\sigma \in
	(0,\sigma _{0})$, we choose 
	\begin{eqnarray}
		&&\lambda =\lambda \left( \sigma \right) =\ln \sigma ^{-1/\left( 2\left(
			2\alpha b-\mu \right) \right) }>\lambda _{2},  \label{3.24} \\
		&&\gamma =\gamma \left( \sigma \right) =2e^{-\lambda \alpha \widetilde{T}%
		}=2\sigma ^{\left( \alpha \widetilde{T}\right) /\left( 2\left( 2\alpha b-\mu
			\right) \right) }.  \label{3.25}
	\end{eqnarray}%
	and let the regularization parameter $\gamma \in \lbrack 2e^{-\lambda \alpha 
		\widetilde{T}},1).$ Let $V_{\min ,\lambda ,\gamma }\left( Y,\tau \right) $
	be the unique minimizer of the functional $J_{\lambda ,\gamma }\left(
	V\right) $ on the set $B\left( r,q_{0},q_{1}\right) $. Let the function $%
	p_{\min ,\lambda ,\gamma }\left( Y\right) $ be defined as: $p_{\min ,\lambda
		,\gamma }\left( Y\right) =4\partial _{Y}V_{\min ,\lambda ,\gamma }\left(
	Y,0\right) $, see (\ref{pp2}). Then there exists a constant $C=C(\mathcal{K}%
	,\alpha ,\mu ,b,r)>0$ depending only on listed parameters such that 
	\begin{eqnarray}
		&&\left\Vert V_{\min ,\lambda ,\gamma }-V^{\ast }\right\Vert _{H^{1}\left( 
			\mathcal{K}_{\alpha ,\mu ,b}\right) }\leq C\sigma ^{\rho },  \label{3.26} \\
		&&\left\Vert p_{\min ,\lambda ,\gamma }-p^{\ast }\right\Vert _{L^{2}\left(
			0,2\alpha b-\mu \right) }\leq C\sigma ^{\rho }.  \label{3.27}
	\end{eqnarray}
\end{theorem}

Now, we state a theorem about the global convergence of the gradient descent
method of the minimization of the functional $J_{\lambda ,\gamma }\left(
V\right) $. Let $\omega \in (0,1)$ be the step size of this method. Fix an
arbitrary number $\vartheta \in (0,1/3)$. We restrict the starting point,
denoted by $V^{(0)}$, to be an arbitrary point in the set $B\left( \vartheta
r,q_{0},q_{1}\right) \subset B(r,q_{0},q_{1})$. Then the gradient descent
method reads as: 
\begin{equation}
	V^{\left( n\right) }=V^{\left( n-1\right) }-\omega J_{\lambda ,\gamma
	}^{\prime }\left( V^{\left( n-1\right) }\right) ,\quad n=1,2,\ldots ,
	\label{gradient}
\end{equation}%
where $V^{(n)}$ approaches the minimizer $V_{\min ,\lambda ,\gamma }$ as $n$
is large. This scheme is well-defined because $J_{\lambda ,\gamma }^{\prime
}\left( V^{(n-1)}\right) \in H_{0}^{4}(\mathcal{K})$; see Theorem \ref%
{sconvex}. Hence, functions $V^{(n)}$ satisfy the same boundary conditions
as in \eqref{ball} for all $n$.

Below, we assume that our minimizer $V_{\min ,\lambda ,\gamma }$ belongs to
the set $B\left( \vartheta r,q_{0},q_{1}\right) $, which is an interior
point of the set $B\left( r,q_{0},q_{1}\right) $. This enabled us in \cite%
{Klibanov2021a} to prove the global convergence of the gradient descent
scheme (\ref{gradient}). This assumption is plausible because by (\ref{3.26}%
) the distance between the ideal solution $V^{\ast }$ and the minimizer $%
V_{\min ,\lambda ,\gamma }$ is small with respect to the noise level $\sigma
\in (0,1),$ although only in the $H^{1}\left( \mathcal{K}_{\alpha ,\mu
	,b}\right) -$norm rather than in the stronger norm $H^{4}\left( \mathcal{K}%
\right) $, see (\ref{ball}).

\begin{theorem}[Global convergence of the gradient descent method 
	\protect\cite{Klibanov2021a}]
	\label{thm:11} Suppose that the parameters $\lambda $ and $\gamma $ are
	taken as in Theorem \ref{sconvex}. Then there exists a sufficiently small
	constant $w_{0}\in (0,1)$ such that the sequence of the gradient descent
	method $\left\{ V^{(n)}\right\} _{n=0}^{\infty }\subset B\left(
	r,q_{0},q_{1}\right) $ for every $\omega \in (0,\omega _{0})$. Moreover, for
	every $\omega \in (0,\omega _{0})$ there exists a number $\theta =\theta
	(\omega )\in (0,1)$ such that 
	\begin{equation*}
		\left\Vert V_{\min ,\lambda ,\gamma }-V^{\left( n\right) }\right\Vert
		_{H^{4}\left( \mathcal{K}\right) }\leq \theta ^{n}\left\Vert V_{\min
			,\lambda ,\gamma }-V^{\left( 0\right) }\right\Vert _{H^{4}\left( \mathcal{K}%
			\right) },\quad n=1,2,\ldots
	\end{equation*}
\end{theorem}

Combining Theorems \ref{thm:10} and \ref{thm:11}, we can prove that the
functions $V^{(n)}$ converge to the ideal solution $V^{\ast }$ as long as
the noise level $\sigma $ in the data tends to zero. In addition, we obtain
the convergence of the corresponding sequence $\left\{ p^{(n)}\right\}
_{n=0}^{\infty }$ towards the ideal function $p^{\ast }$. It works via the
expression $p^{(n)}(Y)=4\partial _{Y}V_{Y}^{(n)}(Y,0)$ in (\ref{pp2}).

\begin{theorem}[\protect\cite{Klibanov2021a}]
	\label{thm:bigthm} Suppose all conditions of Theorems \ref{thm:10} and \ref%
	{thm:11} hold. Then the following convergence estimate is valid: 
	\begin{align*}
		\left\Vert V^{\ast }-V^{\left( n\right) }\right\Vert _{H^{1}\left( \mathcal{K%
			}_{\alpha ,\mu ,b}\right) }& +\left\Vert p^{\ast }-p^{\left( n\right)
		}\right\Vert _{L^{2}\left( 0,2\alpha b-\mu \right) } \\
		& \leq C\sigma ^{\rho }+C\theta ^{n}\left\Vert V_{\min ,\lambda ,\gamma
		}-V^{\left( 0\right) }\right\Vert _{H^{4}\left( \mathcal{K}\right) },\quad
		n=1,2,\ldots ,
	\end{align*}%
	where the constant $C=C(\mathcal{K},\alpha ,\mu ,b,r)>0$ depends only on
	listed parameters.
\end{theorem}

It is now worth mentioning that in our convergence results above, $r$ is arbitrary and the starting
point $V^{(0)}$ is arbitrary in $B(\vartheta r,q_{0},q_{1})$ with $\vartheta 
$ being fixed in $(0,1/3)$. Therefore, our numerical method for solving the underlying nonlinear inverse problem does not require any advanced knowledge of a small neighborhood of the true solution. In this sense, the method is globally convergent, see Definition 1 in Introduction.

\begin{remark}
	Even though the values of $\lambda$ in our CWF should be large in our theory, our rich computational experience working with the convexification tells us that we can choose $\lambda$ between 1 and 3. For this  We refer the reader to our previous publications on the convexification; see \cite{Klibanov2019,Khoa2020,Khoa2020a,Khoa2020b,Klibanov2021a,Smirnov2020,Klibanov2018,Klibanov2019b} some other references cited therein. In the present work, we choose $\lambda=1.05$, and it works well.
\end{remark}

\section{Numerical studies}

\label{sec:experiments}

To generate the data for the inverse problem, we have solved the forward
problem (\ref{forward1})-(\ref{3000}) by the standard implicit scheme.
Therefore, we describe in this section our numerical solution of only the
inverse problem.

It is worth noting that we suppose to have a buried object in the sand
region $\left\{ z<0\right\} \cap \Omega $ of our computational domain $%
\Omega $, see (\ref{eps}). The dielectric constant of this object is assumed
to be different from that of the dry sand. Hence, the sand region is
actually heterogeneous. The solution of problem (\ref{forward1}) includes
the signal from the sand, which causes a significant challenge when working
on the inversion. In this paper, we follow the heuristic approach commenced
in works of this research group with experimentally collected data \cite%
{Klibanov2019,Khoa2020a} to get rid from the sand signals. To do so, for
each source location, we consider $u_{\text{ref}}(\mathbf{x},t)$ as the
reference data. The reference data are generated for the case when an
inclusion is not present in the simulated sandbox. This means that we take $%
\varepsilon _{r}(\mathbf{x})=4$ for $z\leq 0$ and $\varepsilon _{r}(\mathbf{x%
})=1$ for $z>0.$ Indeed, the dielectric constant of the dry sand equals 4
(cf. \cite{Khoa2020a,Khoa2020b}). We remark that we only need to generate
these reference data once for any numerical examples. Let $u^{\ast }(\mathbf{%
	x},t)$ be the simulated data in the case when an inclusion is present in the
sandbox. Then, our target data after filtering out the sand signals are
defined by 
\begin{equation*}
	u_{\text{sc}}(\mathbf{x},t)=u^{\ast }(\mathbf{x},t)-u_{\text{ref}}(\mathbf{x}%
	,t).
\end{equation*}%
And this is the function we work with to solve our inverse problem.

\subsection{Parameters }

\label{sec:5.1}

In this part, we show numerically how our proposed method works to create
3-D images from SAR data taken on three lines $L_{j}^{\text{scr}},j=1,2,3.$
Experimentally, the data collection is very expensive. Therefore, it is
pertinent to use only three lines in these numerical studies. This prepares
a playground to work with experimental SAR data in the next publication.
Different from \cite{Klibanov2021} in which the Lippmann--Schwinger equation
was used to generate the data, our SAR data are computed by the finite
difference solver of section \ref{sec:5.1}. Cf. \cite%
{Klibanov2021,Klibanov2021a}, we use the standard finite difference
approximation to compute the minimizer of $J_{\lambda ,\gamma }$ defined in (%
\ref{JJ}). In this regard, the discrete operators are formulated to
approximate the operator $\mathcal{S}$. Since the integration in $J_{\lambda
	,\gamma }$ is over the two-dimensional rectangle $\mathcal{K}$, we introduce
its equidistant step size as $\left( \Delta _{Y},\Delta _{\tau }\right) $
for the grid points $\left( Y,\tau \right) $ in that domain. Then the
discrete solution is sought by the minimization of the corresponding
discrete cost functional with respect to the values of the function $V$ at
grid points. As in \cite{Klibanov2021}, we only consider the regularization
term in the $H^{2}(\mathcal{K})-$norm for the numerical inversion; compared
with the $H^{4}(\mathcal{K})-$norm in the theory. This reduces the
complexity of computations and thus, saves the elapsed time.

\begin{remark}
	Below we work only with dimensionless variables,
	also, see (\ref{4000}). In this context, 1 in the dimensionless
	spatial variable is 0.3 (meter) in reality, and 1 in the dimensionless time variable is 1 (nanosecond).
\end{remark}

As to our forward problem (\ref{forward1})-(\ref{3000}), we consider $\Omega
=(-2,2)^{3}$. Essentially, our domain of interest is $(-0.6\text{ m},0.6\text{ m})^3$ that corresponds to a scanning ground region of $\left|0.6-(-0.6)\right|^3/2 = 0.864$ ($\text{m}^3$). For $L_{j}^{\text{src}}$, we take $j\in \overline{1,j_{0}}$
with $j_{0}=3$ lines of sources with the increment 0.05 (1.5 cm), which is the
distance between 3 lines of source. The central source position of our
entire configuration is located at $x^{0}=0,y_{2}^{0}=-1.7,z^{0}=0.8$. With
the increment 0.05, our source $y$-locations are $%
y_{1}^{0}=-1.75,y_{2}^{0}=-1.7,y_{3}^{0}=-1.65$, respectively. The length of
each of these lines equals to 1.4 (0.42 \text{m}), i.e. $L=0.7$, and every source line has
28 source positions, i.e. $M=28$. Our circular disk antenna is taken with
the diameter $D=0.1$ (3 cm) and its associated parameters are chosen by $\tau
_{0}=1\;(\text{ns}),\omega _{0}=6\pi \times 10^{8}\;(\text{Hz}),\alpha _{0}=8\pi \times
10^{6}/\tau _{0}$ and $\widetilde{D}=h_{\mathbf{x}}=\left\vert \Omega
\right\vert /(N_{\mathbf{x}}-1)$ with $N_{\mathbf{x}}=81$. As a result, the
carrier wavelength is $\widetilde{\lambda }_{0}=2\pi c_{0}/\omega _{0}=1\;(\text{m})$,
and the Fresnel number $\text{Fr}\approx D^{2}/(\widetilde{\lambda }%
_{0}\left\vert z^{0}\right\vert )=1.25\times 10^{-3}\ll 1$ guarantees well
the far field zone; see Fraunhofer condition in \cite[section 7.1.3]{Lipson2018}. The refinement $N_{\mathbf{x}}$ is chosen because $h_{%
	\mathbf{x}}=0.05$ is compatible with an increment of 1.5 (cm).

\begin{figure}[tbp]
	\begin{centering}
		\subfloat[\label{fig:4a}]{\includegraphics[scale=0.3]{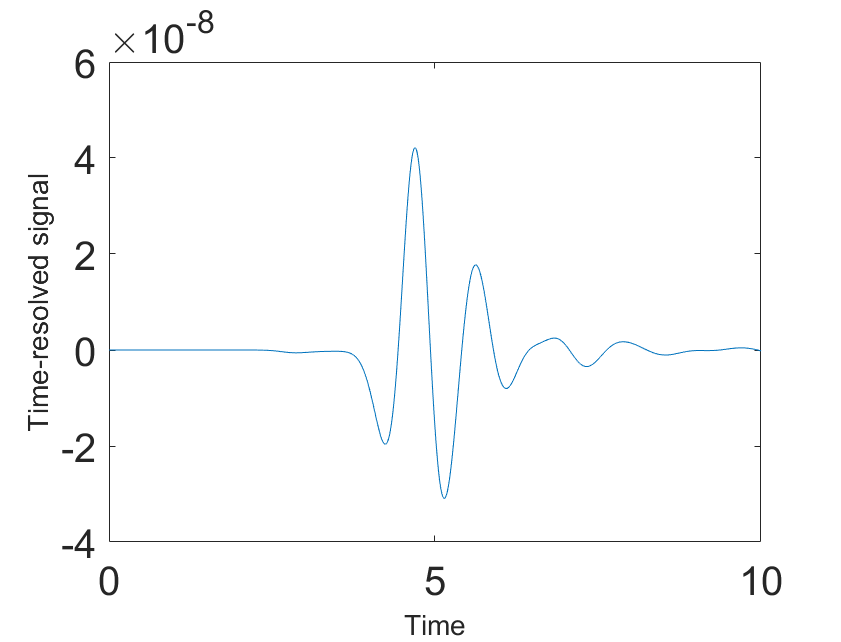}}	
		\subfloat[\label{fig:4b}]{\includegraphics[scale=0.3]{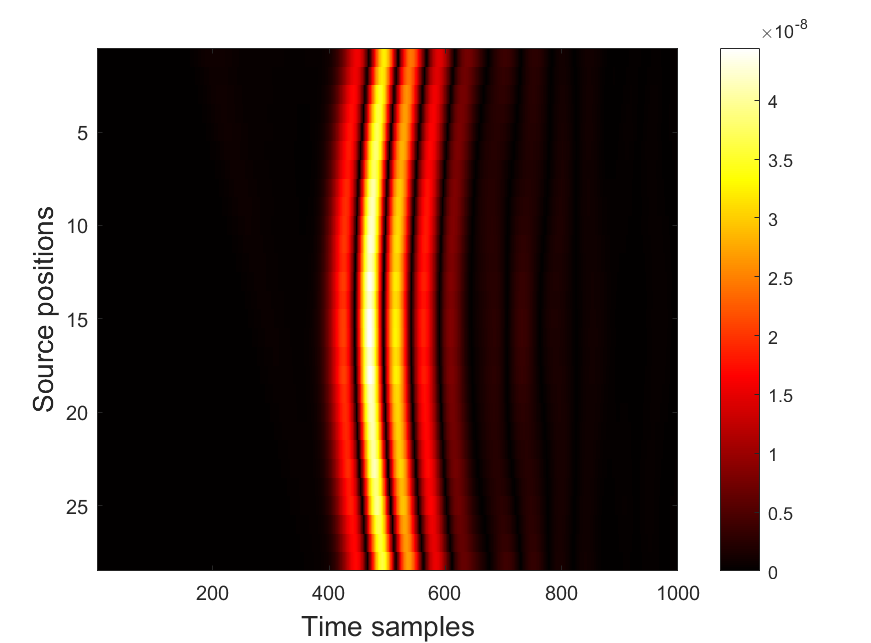}
			
		}
		\par\end{centering}
	\begin{centering}
		\subfloat[\label{fig:4c}]{\includegraphics[scale=0.3]{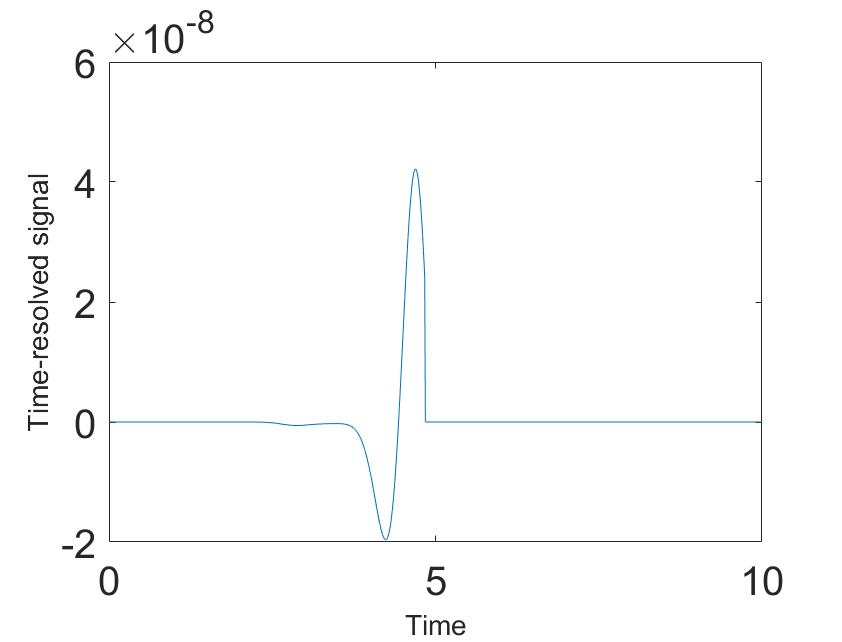}
		}
		\subfloat[\label{fig:4d}]{\includegraphics[scale=0.3]{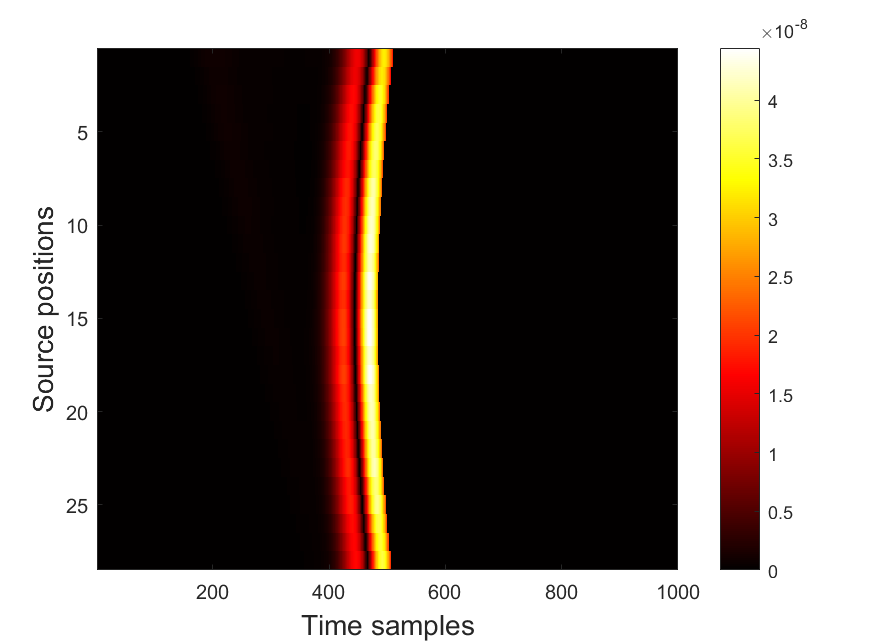}
		}
		\par\end{centering}
	\caption{A typical sample of the time-dependent computationally simulated data for our
		Model 1. This is dimensionless time $\tau = t/\text{ns}$, where $t$ is time in nanoseconds (ns), see (\ref{4000}). (a) The simulated data at the central source position $x^0 = 0,
		y^0_2 = -1.7$ before truncation with $T = 10$. (b) The corresponding
		simulated data with all source positions and 1000 time samples for the time
		domain. (c) The simulated data after truncation. (d) The corresponding
		truncated data with all source positions and 1000 time samples. }
	\label{fig:signal}
\end{figure}

For our numerical results, the parameters were chosen as 
\begin{equation*}
	\lambda =1.05,\alpha =0.49,\gamma =10^{-10}
\end{equation*}
and the discretization $\left( \Delta _{Y},\Delta _{\tau }\right) =(8\times
10^{-3},8\times 10^{-3})$. As mentioned above, even though our theory is valid for large $%
\lambda $, we have observed numerically in all our previous publications on
the convexification that $\lambda \in \lbrack 1,3]$ works well; see our
previous numerical results in, e.g., \cite{Klibanov2019,Khoa2020,Khoa2020a,Khoa2020b,Klibanov2021a,Smirnov2020} and 
\cite{Klibanov2018,Klibanov2019b}. The elevation angle is taken by $\theta
=\pi /4$. Similar to our first work on SAR \cite{Klibanov2021}, we 
multiply the simulated delay-and-sum data by a calibration factor $CF$. Here, we choose $%
CF=1.75\times 10^{15}$. Since our method relies on solving 1D problems, which is an approximate model, we then have no choice but to use that calibration factor.




\subsection{Reconstruction results}

\label{sec:5.2}

We consider two models for our numerical testing. In both cases the objects
are centered at the point $(0,0,-0.14)$.

\begin{itemize}
	\item \textbf{Model 1}: An ellipsoid with principal semi-major axis and two
	semi-minor axes, respectively, being 0.2 (6 cm), 0.12 (3.6 cm) and 0.04 (1.2 cm) in $x,y$ and $z$
	directions respectively. This ellipsoid has the dielectric constant $%
	\varepsilon _{\text{true}}\left( \text{ellip}\right) =15$.
	
	\item \textbf{Model 2}: A rectangular prism with dimensions $0.18\times
	0.08\times 0.08$ (L$\times $W$\times $H) corresponding to $5.4\times 2.4\times 2.4$ in $\text{cm}^3$. It has the dielectric constant $%
	\varepsilon _{\text{true}}\left( \text{prism}\right) =23.8$.
\end{itemize}

\begin{remark}
	We recall in the context of (\ref{trans}) that 1 in the dimensionless
	spatial variable is 0.3 (m) in reality. Therefore, the top point of the
	surface of our simulated objects is close to the air/ground interface $%
	\left\{ z=0\right\} $ within $\left\vert -0.14+0.04\right\vert \times 30=3$
	(cm) of the burial depth, which is typical in de-mining operations; cf. e.g. 
	\cite{Daniels2006}.
\end{remark}

These two examples model well the so-called \textquotedblleft
appearing\textquotedblright\ dielectric constant of metallic objects we
experimented with in \cite{Kuzhuget2012}. It was established numerically in 
\cite{Kuzhuget2012} that the range of the appearing\ dielectric constant of
metals is $[10,30]$. Moreover, this range includes the dielectric constant
23.8 of the experimental watered glass bottle in \cite{Klibanov2021a}. Both
metallic and glass bodied land mines are frequently used in the battlefields.

Suppose that our SAR data are $G(\mathbf{x}_{j}^{0,m},t)$ for each $j=%
\overline{1,j_{0}}$ and $m=\overline{1,N}$. Recall that $\mathbf{x}%
_{j}^{0,m}=\left( x_{m}^{0},y_{j}^{0},z^{0}\right) \in L_{j}^{\text{scr}}$
is the source location number $m$ on the line of sources $L_{j}^{\text{scr}}$
number $j$. Prior of getting the data $g_{0}$ and $g_{1}$ in \eqref{3.4}, we
perform a data preprocessing procedure to filter out unwanted artifacts
observed in the time-resolved signals. First, we eliminate small
\textquotedblleft peaks\textquotedblright\ of the signals by setting: 
\begin{equation}
	\widetilde{G}\left( \mathbf{x}_{j}^{0,m},t\right) =%
	\begin{cases}
		0 & \text{if }\left|G\left( \mathbf{x}_{j}^{0,m},t\right)\right| <0.1\max_{t\in \left[ 0,T%
			\right] ,m\in \lbrack 1,M]}\left\vert G\left( \mathbf{x}_{j}^{0,m},t\right)
		\right\vert , \\ 
		G\left( \mathbf{x}_{j}^{0,m},t\right)  & \text{otherwise}.%
	\end{cases}
	\label{Gtilde}
\end{equation}%
Second, we use the built-in function \texttt{findpeaks} in MATLAB to count
the number of the peaks. Then, we deliberately keep only two first peaks of
the signals because we believe that those are the most important ones for
the object. The resulting signals are denoted by $\widehat{G}(\mathbf{x}%
_{j}^{0,m},t)$. We continue our preprocessing procedure by applying the
delay-and-sum technique to $\widehat{G}(\mathbf{x}_{j}^{0,m},t)$ as in \cite%
{Klibanov2021}. For brevity, we do not detail internal steps of this data
preprocessing procedure, but refer to \cite[section III-A]{Klibanov2021}.

Now we focus on the preprocessing on the delay-and-sum data. In this
context, for each line $L_{j}^{\text{scr}},$ $j=1,2,3$ of source, the data
are along the slant range planes which we reconstruct. Let $P_{j}$ be the slant range plane
which corresponds to the line of sources $L_{j}^{\text{scr}},j=1,2,3.$ In this plane $P_j$, we, slightly abusing the same notations, consider $x_s$ as the source number and $y_s$ as the variable for the axis orthogonal to the line of source $L_j^{\text{scr}}$. Those $y_s$ are actually the transformation of $t$ using time of arrival. Thereby, we denote the
delay-and-sum data by $\text{SAD}_{j}(x_{s},y_{s})$.  We now introduce a truncation to control the size of
the computed object. The \textquotedblleft heuristically\textquotedblright\
truncated delay-and-sum data are given by 
\begin{equation}
	\text{SAD}_{j}^{\text{tr}}(x_{s},y_{s})=%
	\begin{cases}
		\text{SAD}_{j}(x_{s},y_{s}), & \text{if }\text{SAD}_{j}(x_{s},y_{s})\geq
		0.95\max_{x_{s},y_{s}}\left\vert \text{SAD}_{j}(x_{s},y_{s})\right\vert , \\ 
		0, & \text{otherwise}.%
	\end{cases}
	\label{trunc}
\end{equation}%
This way, we use only the absolute values as the data $g_{0}$ for each 1D
inverse problem along the central line of the antenna. We also remark that
the argument of \cite[Remark 4]{Klibanov2021a} implies that the Neumann data 
$g_{1}$ equals to the derivative of $g_{0}$, i.e. $g_{1}=g_{0}^{\prime }.$
Our data preprocessing is exemplified in Figures \ref{fig:4a}--\ref{fig:4d}.

To illustrate how our forward model works, consider Model 1 as an example.
This object is centered at $(0,0,-0.14)$. The central point of the source
line $y_{2}^{0}=-1.7$ is located at $(0,-1.7,0.8)$. Therefore, the distance
between these two points is $\sqrt{(-1.7)^{2}+(0.8+0.14)^{2}}\approx 1.94$.
Hence, the backscattering time-resolved signal should bump at the point $%
1.94\times 2=3.88$ of the time domain. This is approximately the point we
observe in Figure \ref{fig:4a}. More precisely, our approximate signal bumps
at the point $\tau =4.13$, which is 10\% difference with $\tau =3.88$. 
But given that we have the whole ellipsoid rather than just its
central point, this difference is acceptable for our simulations. This
confirms how reliable our forward solver is. We can stop looking for the
wave after the object's bumps, which is after around the point $\tau =5$;
see again Figure \ref{fig:4a}. However, we enlarge the time domain up to the
point $\tau =T=10$ to avoid any possible boundary reflections of the wave
propagation. Moreover, we take into account 1000 time samples in Figure \ref%
{fig:4b} to identify well the location of every wave with respect to $t$.

For each line $L_{j}^{\text{scr}},$ $j=1,2,3$ of sources, we use the
truncation (\ref{trunc}) to form a rectangular area $\left( x_s,y_s\right) \in
\lbrack s_{1}^{j},s_{2}^{j}]\times \lbrack l_{1}^{j},l_{2}^{j}]\subset P_{j}$. As mentioned above when considering variables $x_s,y_s$, this rectangle is used in the slant range plane involving the transformation of $\tau $ in $y$
using times of arrival. 

\begin{remark}\label{rem:5}
	The selection of this rectangular region explains
	the reason why our computed objects in Figures \ref{fig:e2}, \ref{fig:e3}
	and Figures \ref{fig:p2}, \ref{fig:p3} have the rectangle-like shape. We are
	not focusing our work on the shape of the object since it is very
	challenging. Rather, we are interested in the accuracy of the computed
	dielectric constant $\widetilde{\varepsilon }_{r}$ of our two
	models. We are also interested in the dimensions of the computed object, and
	the presence of the above-mentioned rectangular region is helpful in
	controlling those dimensions. 
\end{remark}

Given a slant range plane $P_{j},$ we solve $N$ 1-D CIPs for $N$ sources $%
\mathbf{x}_{j}^{0,m}=\left( x_{m}^{0},y_{j}^{0},z^{0}\right) \in L_{j}^{%
	\text{scr}},m=1,...,N.$  Let the corresponding discrete solutions be
discrete functions defined in  $P_{j},$ $\widetilde{\varepsilon }%
_{1,j}(x_{m}^{0},y_s)$. Note again that slightly abusing notation, the variable $y_s$ here is the grid point in $P_j$. To reduce the artifacts, we average all the 1D
solutions in our rectangular area $\lbrack s_{1}^{j},s_{2}^{j}]\times \lbrack l_{1}^{j},l_{2}^{j}]\subset P_{j}$. To do so, for each $j$, we define first the semi-discrete
function  $\widetilde{\varepsilon }_{2,j}(x_{m}^{0},y)$ for $x_m^0\in[s_1^j,s_2^j],y\in [l_1^j,l_2^j]$ as 
\begin{equation}
	\widetilde{\varepsilon }_{2,j}\left( x_{m}^{0},y\right) =%
	\begin{cases}
		30, & \text{if }\max_{y_s\in [l_1^j,l_2^j]}\left( \left\vert \widetilde{\varepsilon }%
		_{1}\left( x_{m}^{0},y_s\right) \right\vert \right) >500, \\ 
		\max_{y_s\in [l_1^j,l_2^j]}\left( \left\vert \widetilde{\varepsilon }_{1}\left(
		x_{m}^{0},y_s\right) \right\vert \right) , & \text{otherwise,}%
	\end{cases}
	\label{A}
\end{equation}%
since our domain of interest of the dielectric constant is [10,30]. Thus, for each source line of sources  $L_{j}^{\text{scr}}$, 
our computed slant-range solution in $P_j$, denoted by $\widetilde{\varepsilon }_{r,j}(x_{m}^{0},y_s)$,
is computed by
\begin{equation}
	\widetilde{\varepsilon }_{r,j}\left( x_{m}^{0},y_s\right) =%
	\begin{cases}
		\left\vert s_{2}-s_{1}\right\vert ^{-1}\sum_{x_{m}^{0}\in \left[ s_{1}^{j},s_{2}^{j}%
			\right] }\widetilde{\varepsilon }_{2,j}\left( x_{m}^{0},y\right)  & \text{in 
		}\left[ s_{1}^{j},s_{2}^{j}\right] \times \left[ l_{1}^{j},l_{2}^{j}\right] , \\ 
		1 & \text{elsewhere}.%
	\end{cases}
	\label{B}
\end{equation}%
We define our computed slant-range solution in each slant range plane $P_j$ by (\ref{B}) because from the above preprocessing data, we expect our inclusion in only the rectangle area. Therefore, it is relevant that the dielectric constant outside of that area should be unity. Finally, having all computed slant-range solutions $\widetilde{\varepsilon }_{r,1},\widetilde{\varepsilon }_{r,2},\widetilde{\varepsilon }_{r,3}$, we assign one value of the final dielectric
constant $\widetilde{\varepsilon }_{r}$ of our object by averaging the
maximum values of those solutions. In particular, let $\widetilde{\varepsilon }_{\max,j} = \max_{(x_m^0,y_s)\in P_j}\widetilde{\varepsilon }_{r,j}(x_m^0,y_s)$ and let $\text{SL}\subset\mathbb{R}^3$ be the prism formed by 3 slant-range planes $P_j$. We then compute
\begin{equation}
	\widetilde{\varepsilon}_{r}\left(x,y,z\right)=\begin{cases}
		\frac{1}{3}\left(\widetilde{\varepsilon}_{\max,1}+\widetilde{\varepsilon}_{\max,2}+\widetilde{\varepsilon}_{\max,3}\right) & \text{in }\Omega\cap\text{SL},\\
		1 & \text{elsewhere}.
	\end{cases}.  \label{C}
\end{equation}
In addition, we do the linear interpolation with respect
to $j$ to help improve the resolution of images. 

\subsection{Computational results}

\label{sec:5.3}

We now first bring in our results for the case of noiseless data and then for
the case of noisy data.

\subsubsection{Noiseless data}

\label{sec:5.3.1}

In the case of Model 1, our computed dielectric constant for the ellipsoid $%
\widetilde{\varepsilon }_{r}(\text{ellip})=15.06$, compared with the true
value $\varepsilon _{\text{true}}(\text{ellip})=15$. For Model 2, the
computed dielectric constant for the rectangular prism is $\widetilde{%
	\varepsilon }_{r}(\text{prism})=22.97$, compared with the true value $%
\varepsilon _{\text{true}}(\text{prism})=23.8$. Therefore, our
approximations of the dielectric constants are quite accurate. 

\begin{figure}[tbp]
	\begin{centering}
		\subfloat[True\label{fig:e1}]{\includegraphics[scale=0.3]{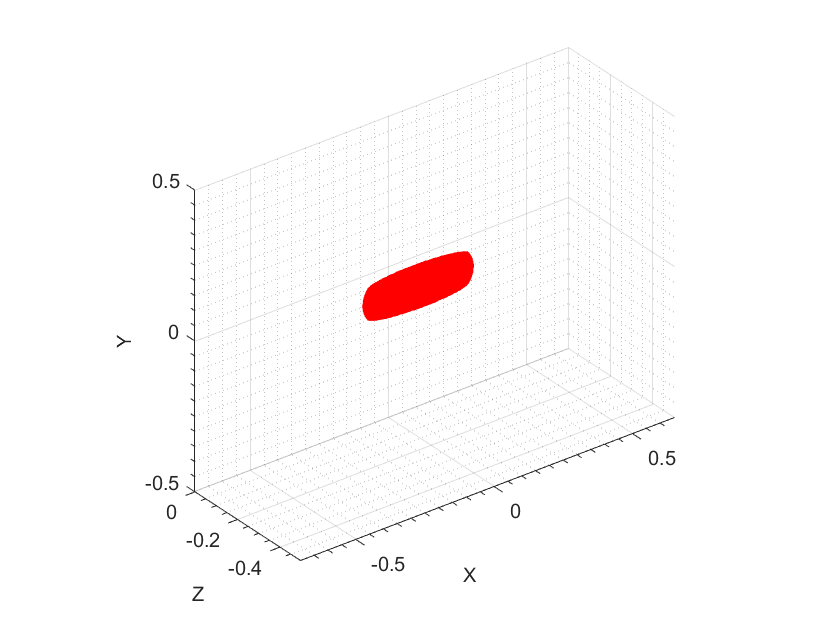}
		}
		\subfloat[Computed\label{fig:e2}]{\includegraphics[scale=0.3]{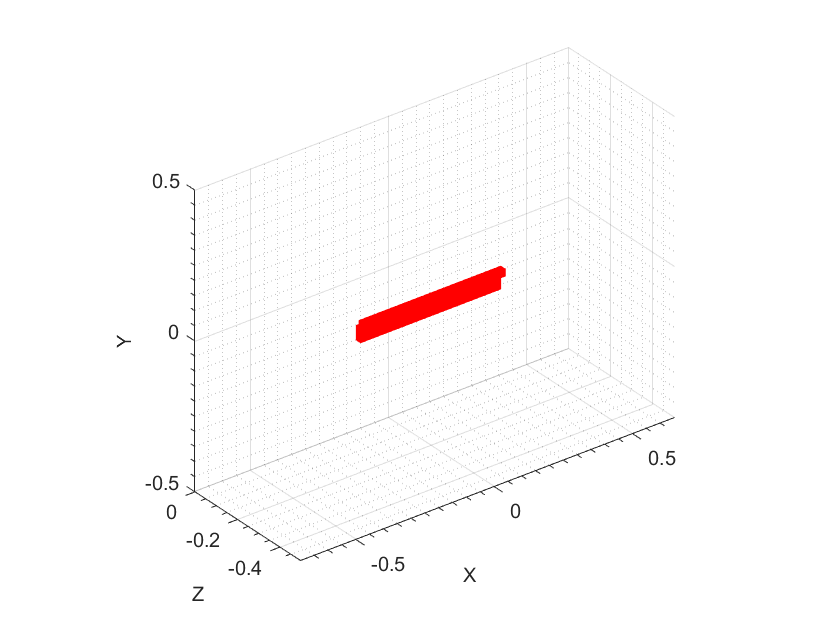}
		}
		\par\end{centering}
	\begin{centering}
		\subfloat[Reconstructed slant-range image\label{fig:e3}]{\includegraphics[scale=0.35]{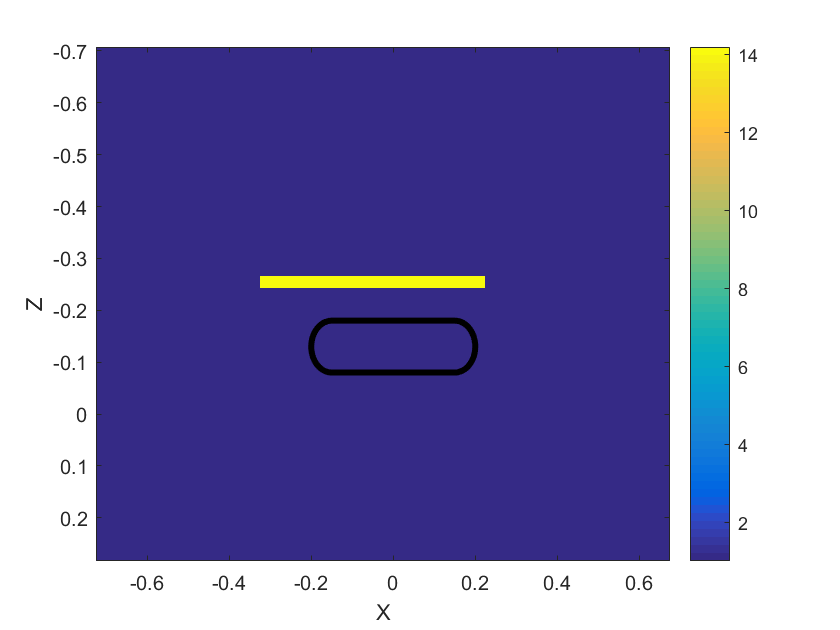}
		}
		\par\end{centering}
	\caption{Our 3-D reconstruction of the ellipsoid of Model 1 with the true
		dielectric constant $\protect\varepsilon_{\text{true}}(\text{ellip})=15$. This is the case of noiseless data. 
		The ellipsoid is buried in a sandbox with the dielectric constant $\protect%
		\varepsilon(\text{sand})=4$. The data are taken in the air with the
		dielectric constant $\protect\varepsilon(\text{air})=1$. (a) True 3D image
		with dimensions $0.4\times0.08\times0.12$ (L$\times$W$\times$H).
		(b) Reconstructed image of the prism with $\protect\widetilde{\protect%
			\varepsilon}_{\text{comp}}(\text{ellip})=15.06$. The \texttt{isovalue} 14 is
		used. Sizes of the computed ellipsoid are $
		0.55\times 0.04\times 0.08$ (L$\times $W$\times $H). Hence, the computed object is longer in length, but smaller in width and is the same in height, as compared with the true ellipsoid.  (c) The cross-section of the 3-D image by the plane $%
		\left\{y=0\right\} $ of the object superimposed with the true one (solid curve). The image with 5\% noise in the data (see (\ref{1000})) is about the same. The computed dielectric constant in the case of noisy data is  $%
		\widetilde{\varepsilon }_{\text{comp}}\left( \text{ellip}\right) =13.52$.}
\end{figure}

\begin{figure}[tbp]
	\begin{centering}
		\subfloat[True\label{fig:p1}]{\includegraphics[scale=0.3]{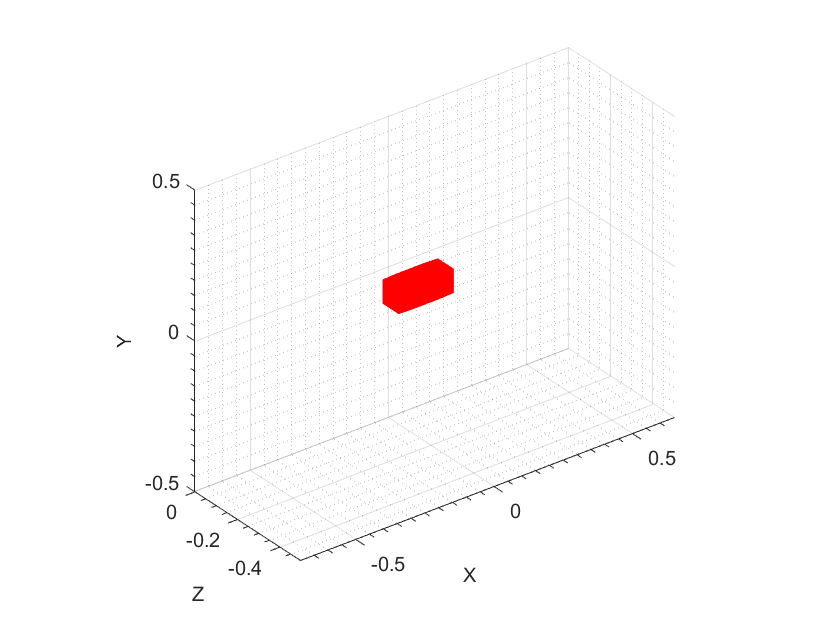}
		}
		\subfloat[Computed\label{fig:p2}]{\includegraphics[scale=0.3]{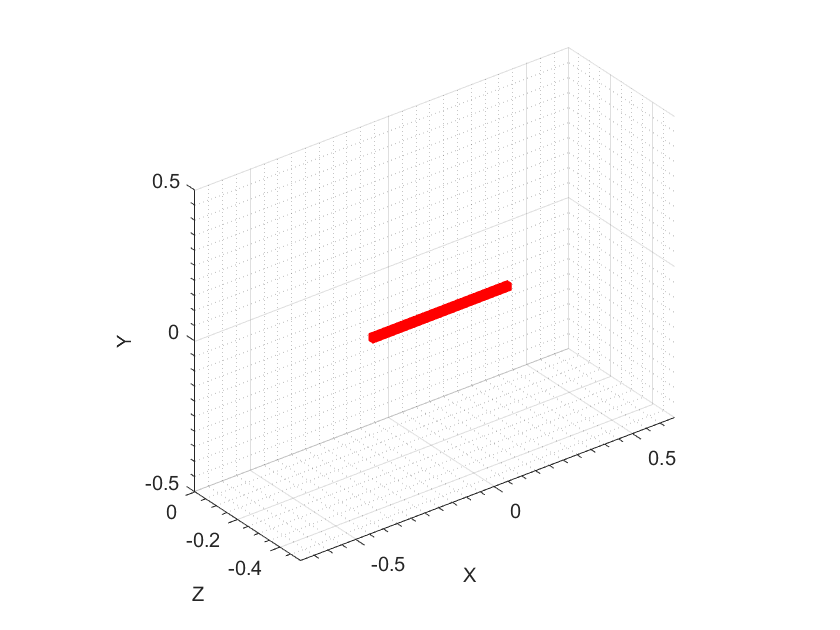}
		}
		\par\end{centering}
	\begin{centering}
		\subfloat[Reconstructed slant-range image\label{fig:p3}]{\includegraphics[scale=0.3]{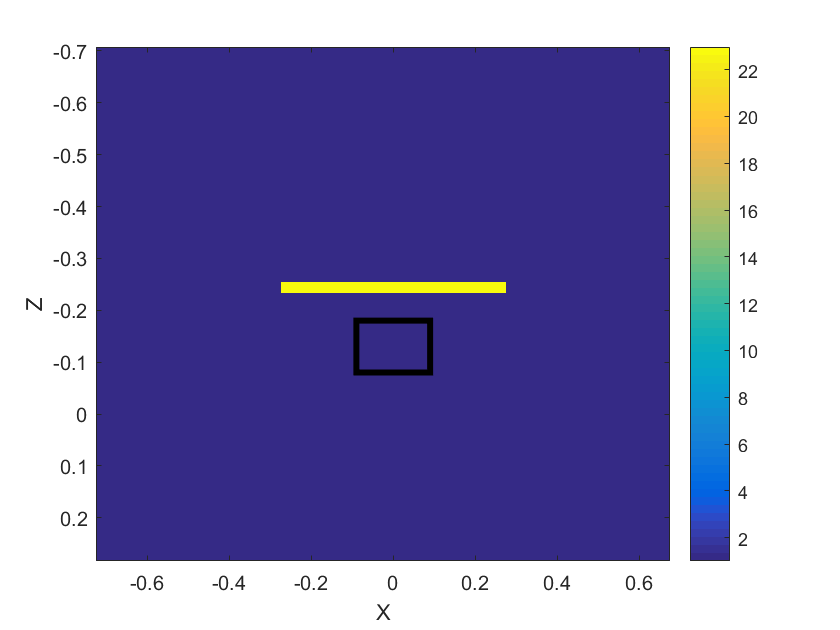}
		}
		\par\end{centering}
	\caption{Our 3-D reconstruction of the rectangular prism of Model 2 with the
		true dielectric constant $\protect\varepsilon_{\text{true}}(\text{prism}%
		)=23.8$. This is the case of noiseless data. The prism is buried in a sandbox with the dielectric constant $%
		\protect\varepsilon(\text{sand})=4 $. The data are taken in the air with the
		dielectric constant $\protect\varepsilon(\text{air})=1$. (a) True 3D image
		with dimensions $0.18\times0.08\times0.08$ (L$\times$W$\times$H). (b)
		Reconstructed image of the prism with $\protect\widetilde{\protect\varepsilon%
		}_{\text{comp}}(\text{prism})=22.97$. The \texttt{isovalue} 22.5 is used. Sizes of the computed prism are $0.5\times 0.03\times 0.03$ (L$\times $W$\times $H). Hence, the computed
		prism is longer in length but smaller in width and height than the true
		prism. 
		(c) The cross-section of the 3-D image by the plane $\left\{y=0\right\}$ of
		the object superimposed with the true one (solid curve). The image with 5\% noise in the data (see (\ref{1000})) is about the
		same. The computed dielectric constant in the case of noisy data is $%
		\widetilde{\varepsilon}_{\text{comp}}\left( \text{prism}\right) =22.08.$}
	\label{fig:jug}
\end{figure}

We now comment on dimensions of computed images. Recall that by
Remark \ref{rem:5} our images have only rectangular shapes. One can see in Figures \ref{fig:e1} and \ref{fig:e2} that the dimensions of
the computed ellipsoid of Model 1 are about $0.55\times 0.04\times 0.08$ ($%
\text{L}\times \text{W}\times \text{H}$). On the other hand, dimensions of
the true ellipsoid are $0.4\times 0.24\times 0.08$ ($\text{L}\times \text{W}%
\times \text{H}$). This means that the computed object is longer in length,
but smaller in width and is the same in height. Similarly, for Model 2, one
can see  in Figures \ref{fig:p1} and \ref{fig:p2} that the dimensions of the
computed prism are $0.5\times 0.03\times 0.03$ ($\text{L}\times \text{W}%
\times \text{H}$), while dimensions of the true prism are $0.36\times
0.16\times 0.16$ ($\text{L}\times \text{W}\times \text{H}$)$.$ Hence, the
computed prism is longer in length but smaller in width and height than the
true prism. 

We remark that the \texttt{isovalue} for the presentation of the dielectric
constants (as well as sizes of imaged targets) in Figures \ref{fig:e2} and \ref{fig:p2} equals to $\varepsilon _{%
	\text{true}}-1$ because we want to see how close our approximate dielectric
constant to the true one is. 

Consider now  Figures \ref{fig:e3} and \ref{fig:p3} in which we have the 2D
cross-sections of the computed images by the plane $\left\{ y=0\right\} $,
the location is not so accurate as we have 10\% difference in location from
our forward solver. We hope to improve these in our future work.

\subsubsection{Noisy data}

\label{sec:5.3.2}

While the above results are for the case of the\ noiseless data, we now
introduce a noise in the data. We add a random additive noise to the
truncated simulated data $\widetilde{G}\left( \mathbf{x}_{j}^{0,m},t\right) $
(see (\ref{Gtilde})). In this regard, we introduce 
\begin{equation}
	\widetilde{G}_{\text{noise}}\left( \mathbf{x}_{j}^{0,m},t\right) =\widetilde{%
		G}\left( \mathbf{x}_{j}^{0,m},t\right) +\sigma \text{rand}\left( \mathbf{x}%
	_{j}^{0,m},t\right) \max_{\mathbf{x}_{j}^{0,m},t}\left\vert \widetilde{G}%
	\left( \mathbf{x}_{j}^{0,m},t\right) \right\vert ,\quad j=1,2,3,m=1,...,N.  \label{1000}
\end{equation}%
Here, $\sigma \in (0,1)$ represents the noise level and \textquotedblleft
rand\textquotedblright\ is a random number uniformly distributed in the
interval $(-1,1)$. We use now the noisy data (\ref{1000}) for both Models 1
and 2 with $\sigma =0.05,$ which is $5\%$ of noise. 

The resulting computed dielectric constant for the ellipsoid is $\widetilde{%
	\varepsilon }_{r}(\text{ellip})=13.52$, compared with the true value $%
\varepsilon _{\text{true}}(\text{ellip})=15$. For Model 2, our computed
dielectric constant for the rectangular prism is $\widetilde{\varepsilon }%
_{r}(\text{prism})=22.08$, compared with the true value $\varepsilon _{\text{%
		true}}(\text{prism})=23.8$. In both models images for noisy data are
about the same as the ones on Figures 3 and 4 for noiseless data.

\section{Concluding remarks}

\label{sec:6}

In summary, we have numerically tested the convexification-based nonlinear
SAR (CONSAR) imaging on identifying images of the dielectric constant of
buried objects in a three-dimensional setting. We have numerically observed
provides accurate values of dielectric constants of targets. However, the
computed locations and sizes of tested targets are not as accurate as the
convexification usually provides for conventional CIPs, see, e.g. the works
of this research group for applications oof the convexification to
experimental data \cite{Khoa2020a,Khoa2020b}. We hope to address these
questions in follow up publications.

\section*{Acknowledgment}

Vo Anh Khoa acknowledges Professor Taufiquar Khan (University of North
Carolina at Charlotte) for the moral support of his research career.

\end{document}